\documentclass[12pt]{amsart}

\usepackage{booktabs}

\makeatletter
\@namedef{subjclassname@2020}{\textup{2020} Mathematics Subject Classification}

\usepackage{youngtab}

\makeatother

\newcommand{\ssyt}{\te{SSYT}}

\newcommand{\row}{\text{row}\,}

\newcommand{\rowt}{\text{row}\,T}
\renewcommand{\lg}{*(lightgray) }

\newcommand{\partition}{\vdash}
\newcommand{\parn}{\vdash n}

\newcommand{\qte}[1]{\q\te{#1}}

\renewcommand{\b}{\mathbf{b}}

\newcommand{\oll}{\ol{\el}}

\newcommand{\R}{\mathcal{R}}

\newcommand{\dsum}{\di\sum}

\newcommand{\dep}[1]{{\vrule width 0pt height 0pt depth #1}}

\newcommand{\el}{\ensuremath{\ell}}

\DeclareFontFamily{U}{mathx}{\hyphenchar\font45}
\DeclareFontShape{U}{mathx}{m}{n}{
      <5> <6> <7> <8> <9> <10>
      <10.95> <12> <14.4> <17.28> <20.74> <24.88>
      mathx10
      }{}
\DeclareSymbolFont{mathx}{U}{mathx}{m}{n}
\DeclareFontSubstitution{U}{mathx}{m}{n}
\DeclareMathAccent{\widecheck}{0}{mathx}{"71}

\newcommand{\ch}{\ensuremath{\widecheck}}

\newcommand{\tw}{\textwidth}

\renewcommand{\kill}[1]{}
\newcommand{\dummy}[1]{\mbox{}}

\makeatletter
\newcommand{\xequal}[2][]{\ext@arrow 0055{\equalfill@}{#1}{#2}}
\def\equalfill@{\arrowfill@\Relbar\Relbar\Relbar}
\makeatother

\newcommand{\mto}{\mapsto}

\newcommand{\ku}{\ensuremath{\emptyset}}

\renewcommand{\k}{\ensuremath{\ol{\mathrm{P}}}}

\newcommand{\vf}{\vfill}

\newcommand{\h}{\hline}

\newcommand{\kyo}[1]{{\eh{#1}}}

\newcommand{\m}{\ensuremath{\infty}}

\renewcommand{\k}[1]{\ensuremath{\left({#1}\right)}}

\newcommand{\ds}{\dots}

\newcommand{\bca}{\begin{cases}}
\newcommand{\eca}{\end{cases}}

\newcommand{\maj}{\textnormal{maj}}

\newcommand{\des}{\textnormal{des}\,}
\newcommand{\A}{\mathcal{A}}

\newcommand{\bpic}{\begin{picture}}\newcommand{\epic}{\end{picture}}

\newcommand{\beda}{\begin{edaenumerate}}
\newcommand{\eeda}{\end{edaenumerate}}

\newcommand{\g}{\ensuremath{\mathbf{g}}}

\newcommand{\cd}{\cdots}

\newcommand{\xto}{\xrightarrow}

\newcommand{\st}{\strut}

\newcommand{\q}{\quad}

\newcommand{\too}{\longrightarrow}

\newcommand{\bq}{\begin{quote}}\newcommand{\eq}{\end{quote}}
\newcommand{\gam}{\gamma}

\renewcommand{\sp}[1]{\ul{\ph{#1}}}
\renewcommand{\sp}[1]{\text{sp}}

\newcommand{\ti}{\times}

\newcommand{\be}{\begin{enumerate}}\newcommand{\ee}{\end{enumerate}}
\newcommand{\bce}{\begin{center}}\newcommand{\ece}{\end{center}}
\newcommand{\bde}{\begin{description}}\newcommand{\ede}{\end{description}}
\newcommand{\bri}{\begin{flushright}}\newcommand{\eri}{\end{flushright}}
\newcommand{\bb}{\begin{block}}\newcommand{\eb}{\end{block}}
\newcommand{\bt}{\begin{thm}}\newcommand{\et}{\end{thm}}
\newcommand{\bpf}{\begin{proof}}\newcommand{\epf}{\end{proof}}
\newcommand{\bex}{\begin{ex}}\newcommand{\eex}{\end{ex}}
\newcommand{\bexr}{\begin{exr}}\newcommand{\eexr}{\end{exr}}
\newcommand{\bft}{\begin{fact}}\newcommand{\eft}{\end{fact}}
\newcommand{\brk}{\begin{rmk}}\newcommand{\erk}{\end{rmk}}
\newcommand{\ba}{\begin{align*}}\newcommand{\ea}{\end{align*}}
\newcommand{\bexe}{\begin{exe}}\newcommand{\eexe}{\end{exe}}

\newcommand{\bit}{\begin{itemize}}\newcommand{\eit}{\end{itemize}}

\newcommand{\bcm}{}
\newcommand{\ol}{\overline}\newcommand{\ul}{\underline}
\newcommand{\hf}{\hfill}

\newcommand{\cc}{\ensuremath{\mathbf{C}}}
\newcommand{\nn}{\ensuremath{\mathbf{N}}}
\newcommand{\qq}{\ensuremath{\mathbf{Q}}}

\newcommand{\zz}{\ensuremath{\mathbf{Z}}}

\newcommand{\bd}{\begin{defn}}\newcommand{\ed}{\end{defn}}
\newcommand{\bp}{\begin{prop}}\newcommand{\ep}{\end{prop}}

\newcommand{\eh}{\emph}\newcommand{\al}{\alpha}
\newcommand{\sub}{\subseteq}
\newcommand{\lam}{\lambda}

\newcommand{\mb}{\mbox}
\newcommand{\te}{\text}\newcommand{\ph}{\phantom}
\newcommand{\wt}{\widetilde}\newcommand{\sm}{\setminus}
\newcommand{\lef}{\left}\newcommand{\ri}{\right}
\renewcommand{\l}{\left}

\newcommand{\then}{\Longrightarrow}

\newcommand{\di}{\displaystyle}\renewcommand{\a}{\ensuremath{\bm{a}}}
\renewcommand{\b}{\ensuremath{\bm{b}}}
\renewcommand{\d}{\ensuremath{\bm{d}}}
\newcommand{\e}{\ensuremath{\bm{e}}}

\newcommand{\f}{\frac}
\newcommand{\x}{\ensuremath{\bm{x}}}

\renewcommand{\b}{\beta}

\renewcommand{\x}{\mathbf{x}}

\newcommand{\F}{\mathcal{F}}
\renewcommand{\d}{\delta}

\newcommand{\Q}{\mathcal{Q}}

\renewcommand{\d}{\delta}

\usepackage{shuffle}
\usepackage[dvipdfmx]{graphicx}
\usepackage[dvipsnames]{xcolor}
\usepackage{float,afterpage}
\usepackage{asymptote,layout}
\usepackage{wrapfig,epic}
\usepackage{amsmath,amssymb,cases,pict2e}

\usepackage{multicol}
\usepackage{tikz}
\tikzset{
    cell/.style={
        anchor=south west,
        draw,
        minimum size=1cm,
    },
}

\usepackage{geometry}
\usepackage{exscale,latexsym,bm}
\usepackage{amssymb,enumerate,amsmath,amsthm,amsfonts}
\usepackage{verbatim,fancybox,wasysym,fancyhdr,type1cm}
\usepackage[frame,all,poly,curve,knot,arrow]{xy}
\usepackage{ytableau}
\usepackage{colortbl}
\usepackage{boxedminipage}
\usepackage{multirow}

\graphicspath{{./figs/}}
\theoremstyle{definition}
\newtheorem{thm}{Theorem}[section]

\newtheorem{lem}[thm]{Lemma}
\newtheorem{prop}[thm]{Proposition}\newtheorem{cor}[thm]{Corollary}

\newtheorem{exr}[thm]{Exercise}

\newtheorem{ex}[thm]{Example}

\newtheorem{defn}[thm]{Definition}\newtheorem{rmk}[thm]{Remark}
\newtheorem{fact}[thm]{Fact}
\newtheorem{block}[thm]{}
\newtheorem*{exe}{Exercise}

\renewcommand{\a}{\alpha}

\renewcommand{\l}{\lam}
\renewcommand{\h}{\hline}
\renewcommand{\arraystretch}{1.5}

\renewcommand{\g}{\gamma}

\renewcommand{\g}{\Gamma}

\renewcommand{\R}{\mathbf{R}}
\usepackage{youngtab}

\renewcommand{\l}{\lambda}
\renewcommand{\wt}{\text{wt}\,}
\renewcommand{\ep}{\varepsilon}
\newcommand{\yd}{\ydiagram}

\renewcommand{\g}{\gamma}

\renewcommand{\m}{\mu}

\renewcommand{\A}{\mathbf{A}}
\newcommand{\B}{\mathbf{B}}

\renewcommand{\st}{\text{st}\,}

\renewcommand{\a}{\bm{a}}

\renewcommand{\a}{\al}
\renewcommand{\B}{\mathbf{B}}

\renewcommand{\a}{\alpha}

\renewcommand{\b}{\beta}

\renewcommand{\B}{\mathbf{B}}

\renewcommand{\m}{\mu}

\renewcommand{\l}{\lambda}

\renewcommand{\B}{\mathbf{B}}
\renewcommand{\bt}{\te{BT}}

\newcommand{\dest}{\te{des}\,T}

\renewcommand{\oll}{\ol{\lam}}

\newcommand{\tref}{\te{ref}}

\newcommand{\Des}{\te{Des}}

\renewcommand{\F}{\mathbf{F}}
\renewcommand{\Q}{\mathbf{Q}}
\renewcommand{\g}{\gamma}

\renewcommand{\a}{\alpha}
\renewcommand{\b}{\beta}

\renewcommand{\ch}{\te{ch}\,}

\newcommand{\dst}{\te{dst}\,}

\renewcommand{\flat}{\te{flat}}

\newcommand{\qy}{\te{QY}}
\newcommand{\qyl}{\ensuremath{\te{QY}(\l)}}

\renewcommand{\st}{\te{ST}}
\newcommand{\ssytl}{\ensuremath{\te{SSYT}(\l)}}
\newcommand{\stl}{\ensuremath{\te{ST}(\l)}}

\newcommand{\csl}{\ensuremath{\te{CS}(\l)}}

\newcommand{\Sk}{\te{Sk}}
\newcommand{\Skl}{\te{Sk}_{\l}}

\newcommand{\comp}{\models}
\newcommand{\compn}{\models n}

\newcommand{\Parn}{\te{Par}(n)}
\renewcommand{\ph}{\Phi}

\renewcommand{\st}{\te{ST}}
\renewcommand{\wt}{\text{wt}\,}
\newcommand{\lparn}{\l\parn}
\newcommand{\Comp}{\te{Comp}}
\renewcommand{\comp}{\te{comp}}
\renewcommand{\f}{\mathbf{f}}
\renewcommand{\A}{\mathbf{A}}

\newcommand{\K}{\mathbf{K}}
\renewcommand{\dep}{\te{dep}}
\renewcommand{\ch}{d}
\renewcommand{\ch}{\te{dep}}

\newcommand{\depw}{\dep\,w}
\renewcommand{\oll}{\ol{\l}}

\newtheorem{prob}[thm]{Problem}

\begin{document}

\title[crystal skeleton polynomials]
{crystal skeleton polynomials with major index, charge and depth}


\author{Masato Kobayashi}

\date{\today}                                       



\address{Masato Kobayashi\\
Department of Engineering\\
Kanagawa University, Rokkaku-bashi, Yokohama, Japan.}
\email{masato210@gmail.com}
%

\maketitle
\begin{abstract}
We introduce a new family of 
polynomials, crystal skeleton polynomials, 
to better understand enumeration of standard Young tableaux, quasi-Yamanouchi tableaux 
and interactions with Gessel's expansion of 
a Schur function, quasi-crystals and crystal skeletons as Maas-Gari\'{e}py introduced in 2023. 
After developing calculus of those polynomials, 
we organize thoughts on major index, 
charge, depth, inversions with RSK  correspondence and a bivariate factorial. 
Also, we revisit the theorem on internal zeros of 
fake degree polynomials by Billey--Konvalinka--Swanson (2020). 
These results altogether improve Gessel's expansion.
\end{abstract}

\tableofcontents

\ytableausetup{centertableaux}

\tableofcontents

\section{Introduction}
Enumeration of \kyo{standard Young tableaux} (SYT) is a classical topic in algebraic 
combinatorics such as hook-length formulas, 
Specht modules, reduced words of permutations, to name a few. The goal of this article is to introduce a new family of 
polynomials, \kyo{crystal skeleton polynomials}, 
as an auxiliary tool to 
better understand such enumeration. 
First, we explain our motivation briefly.

Let $X=(x_{1}, x_{2}, \ds)$ denote a set of 
infinite variables and $g(X)\in \cc[[X]]$. Suppose 
$g(X)$ is homogeneous of degree $n$.
Say $g$ is \kyo{$F$-positive} if 
\[
g(X)=\sum_{\al\compn}c_{\al}F_{\al}(X), 
\q c_{\a}\in \zz_{\ge0}
\]
where $\al$ is a strong composition of $n$ 
and $F_{\a}$ a fundamental quasi-symmetric function. 
It is a basic question to ask when an $F$-positive function is Schur-positive. 
If $g$ is symmetric, then it must be Schur positive by the theory on the ring of symmetric functions;
See \cite{elw, ge2, ossz, ro} on 
how to find its Schur expansion from such $g$. 
The next question is this.
\begin{quote}
When is an $F$-positive function 
a \kyo{single} Schur function?
\end{quote}
Gessel found the following expansion \cite{ge}
\[
s_{\l}(X)=\sum_{T\in \stl}F_{\dest}(X)
\]
where $\stl$ denotes the set of all SYTs of shape $\l$ 
and $\dest$ is the descent composition of $T$ (not   descent set). 
Assaf--Searles \cite{asse} improved it 
a bit in terms of \kyo{quasi-Yamanouchi tableaux} under polynomial setting. Moreover, 
Assaf \cite{as} gave axioms of a \kyo{dual equivalence graph} which characterizes such a sum, roughly speaking. 
The idea of such a graph originally goes back to Haiman \cite{ha}.

Let us now call a (multi)set 
$\te{Sch}^{\te{(multi)}}(\l):=\{\dest\mid T\in \stl\}$ the \kyo{Schur family} for $\l$. 
We aim to understand a larger collection 
$\mathbf{Sch}=
\{\te{Sch}(\l)\mid \te{$\l\in \te{Par}$}\}$ rather globally. 
Our framework fits in with the brand new theory of  a \kyo{crystal skeleton} which 
Maas-Gari\'{e}py \cite{ma} introduced only in 2023.
Afterward, Cain--Malheiro--Rodrigues--Rodrigues (2025) \cite{cmrr} and Brauner--Corteel--Daugherty--Schilling (2025) \cite{bcds} 
follow her project.
%




We assume familiarity with 
basic ideas on Young tableaux; 
for undefined terms, see \cite{sta}. 
Throughout, a fundamental quasi-symmetric  function 
and a descent \kyo{composition} are key concepts  \cite{bcds}. 
Below, $\l$ is always a partition, 
$n=|\l|$ and $\al$ is a strong composition of size $n$. We sometimes use shorthand notation 
$\a=\a_{1}\ds \a_{l}$ to mean 
$\a=(\a_{1}, \ds, \a_{l})$. 

\section{Skeleton polynomials}

\subsection{Preliminaries}

A \kyo{weak (strong) composition} is a finite sequence 
of nonnegative (positive) integers as we often write 
$a=(a_{1}, \ds, a_{m})$. 
The \kyo{size} of $a$ is 
\[
|a|=a_{1}+\cd+a_{m}.
\]
By $\comp$ ($\Comp$) we mean the set of all 
weak (strong) compositions. 

Let $a, b\in \Comp$ and $c\in \comp$. 
Define $b=\flat(c)$ if we obtain $b$ with 
deleting all 0 parts of $c$. 
Say $b$ \kyo{refines} $a$ if 
we can obtain $b$ from $a$ with repeating 
to split one number into two adjacent positive entries.
Denote 
by $\tref(a)$ the set of all 
strong compositions which refine $a$.
For example, if $a=(2, 3)$, then 
\[
\tref{(a)}=
\{23, 221, 212, 2111, 113, 1121, 1112, 11111\}.
\]

Let 
$X=(x_{1}, x_{2}, \ds)$ be infinitely many variables. 
To any composition $a$, 
associate the weight monomial 
$x^{a}=x_{1}^{a_{1}}x_{2}^{a_{2}}\cd$.


\begin{defn}
Let $a, b\in {\Comp}$. 
The \kyo{monomial quasi-symmetric function} for $b$ is 
\[
M_{b}(\x)=\sum_{
\substack{c\in \comp\\\flat(c)=b}
}
x^{c}.
\]
The \kyo{fundamental quasi-symmetric function} for $a$ is 
\[
F_{a}(X)=\sum_{b\in \tref{(a)}}
M_{b}(X).
\]
\end{defn}
To describe 
connections with such functions and SSYTs, we need another definition; a skew shape is a \kyo{horizontal strip} if it contains at most one box in each column.

\begin{rmk}
In this article, a \kyo{box} always means \kyo{only a position} in a diagram. 
\end{rmk}

Let $B_{1}, B_{2}$ be boxes in a diagram. 
Write $B_{1}<_{nE} B_{2}$ if $B_{2}$ lies 
northEast to $B_{1}$, that is, 
if 
$B_{1}=(i_{1}, j_{1})$ and  
$B_{2}=(i_{2}, j_{2})$, 
then $i_{1}\ge i_{2}$ and $j_{1}<j_{2}$.

Suppose now 
$B=(B_{1}, \ds, B_{k})$ 
is a sequence of boxes and 
moreover there is an associated sequence 
$a=(a_{1}, \ds, a_{k})$ such that $a_{j}\in B_{j}$.



\begin{defn}
Say $H=(B, a)$ is a \kyo{horizontal band} if 
\begin{itemize}
	\item $B, u$ are nonemtpy,
	\item $B_{1}<_{nE}\cd <_{nE}B_{k}$,
	\item $a_{j}\le a_{j+1}$ for all $j$.
\end{itemize}
For convenience, let $|H|=k$.
\end{defn}
%
For an SSYT $T$, 
there exists a unique sequence of maximal horizontal bands $H_{1}, \ds, H_{l}$ 
such that 
they contain each entry of $T$ exactly 
once and if all entries in $H_{i+1}$ are strictly greater than than those in $H_{i}$; 
maximal means adding a box to $H_{i}$ 
cannot be any other horizontal band in $T$. 
We call such a sequence the (semistandard) \kyo{minimal parsing} of $T$. 

\begin{defn}
The \kyo{descent composition} of an SSYT $T$ 
is the strong composition 
\[
\des\, T =\k{|H_{1}|, |H_{2}|, \ds, |H_{l}|}
\]
as $H_{i}$'s described just above. 
%
\end{defn}


Given an SSYT $T$ with 
its descent composition $\al=(\a_{1}, \ds, \al_{l})$
, its minimal parsing allows us to define the 
\kyo{standardization} $T^{st}$:
all entries in the first maximal horizontal band 
are relabeled by $1, \ds, \a_{1}$ (from left), 
in the second, by $\a_{1}+1, \ds, \a_{1}+\a_{2}$ and so on. \kyo{Destandardization} of $T$ 
$(T^{dst})$ is the tableau with replacing all numbers in $H_{i}$ by $i$.

\begin{ex}
Let 
\ytableausetup{mathmode, boxsize=14pt}
$T=
\begin{ytableau}
	1&1&2&5\\
	3&8\\8
\end{ytableau}$. 
Its minimal parsing is 
\[
\k
{\,\begin{ytableau}
	1&1&2&\lg{}\\
	\lg{}&\lg{}\\\lg{}
\end{ytableau}, 
\begin{ytableau}
	\lg{}&\lg{}&\lg{}&5\\
	3&\lg{}\\\lg{}
\end{ytableau}, 
\begin{ytableau}
	\lg{}&\lg{}&\lg{}&\lg{}\\
	\lg{}&8\\8
\end{ytableau}\,}
\]
and $\dest=(3, 2, 2)$, 
$T^{st}=
\begin{ytableau}
	1&2&3&5\\
	4&7\\6
\end{ytableau}$, 
$T^{dst}=
\begin{ytableau}
	1&1&1&2\\
	2&3\\3
\end{ytableau}.$

\end{ex}




\subsection{Quasi-crystal, crystal skeleton}

%
%
%


Our discussion proceeds with $\B(\l)=(\ssyt(\l)\sqcup\{0\}, \e_{i}, \f_{i}, \to)$, a crystal graph of type A. Formally, set $0^{st}=0$.

\begin{defn}
Say $T, U\in \B(\l)$ are 
\kyo{standard equivalent} $(T\sim U)$ if 
$T^{st}=U^{st}$. 
A \kyo{quasi-crystal} is an equivalent class of 
this equivalent relation on $\B(\l)$ excluding $\{0\}$. 
A \kyo{crystal skeleton} $\csl$ is 
the quotient set 
$\B(\l)/\!\!\sim$.
\end{defn}

\begin{rmk}
For here, treat $\csl$ just as a set. We can however impose a graph structure with certain vertex- or edge-labeling and direction of edges \cite{bcds, cmrr, ma}. 
It contains a dual equivalence graph as a connected  subgraph anyway.
%
\end{rmk}
\begin{fact}[\cite{ma}]
Each quasi-crystal in $\B(\l)$ forms a connected subgraph. Its generating function is 
single $F_{\a}(X)$ for some $\a\models |\l|$. 
\end{fact}

\begin{ex}
\ytableausetup{mathmode, boxsize=8pt}
Figure \ref{f1} illustrates 
five quasi-crystals in $\B_{3}\k{\yd{3,2}}$. 
A dotted edge connects distinct quasi-crystals.
\end{ex}
Let us now clarify a deeper 
interpretation of Gessel expansion \cite{ma} with little more new words. If $T\sim U$ in $\B(\l)$, then 
$\des T=\des U$ as we can check. 
Therefore, the descent composition for a quasi-crystal is well-defined.
Say quasi-crystals $\qq, \mathbf{R}\sub\B(\l)$ are \kyo{descent equivalent} if 
$\des \qq=\des \mathbf{R}$. 
In fact, descent equivalent quasi-crystals have the identical generating function.
\begin{defn}
A \kyo{fundamental system}
for $\a$ in $\B(\l)$ is 
\[
\F(\l, \a)=\bigsqcup \Q
\]
where $\Q$ runs all quasi-crystals in $\B(\l)$ such that $\des \Q=\al$.
\end{defn}
All connected components in 
$\F(\l, \a)$ are isomorphic as colored directed graphs.
As a consequence, $\B(\l)$ 
decomposes into a union of such $\F(\l, \a)$'s. 




\begin{figure}
\begin{center}
\ytableausetup{mathmode, boxsize=8pt}
\caption{Quasi-crystals in $\B_{3}\k{\yd{3,2}}$}
\label{f1}
\begin{center}
\mb{}\vf
\vf
\ytableausetup{mathmode, boxsize=15pt}
\begin{minipage}[c]{.85\tw}
\xymatrix@=9mm{
*{}&*+{\begin{ytableau}
	1&1&1\\
	2&2
\end{ytableau}}\ar@{.}_{}[ld]
\ar@{->}_{2}[rd]
&\\
*+{\begin{ytableau}
	1&1&2\\
	2&2
\end{ytableau}}\ar@{->}_{2}[d]
&&*+{\begin{ytableau}
	1&1&1\\
	2&3
\end{ytableau}}\ar@{.}_{}[ld]\ar@{->}_{2}[d]\\
*+{\begin{ytableau}
	1&1&3\\
	2&2
\end{ytableau}}\ar@{->}_{2}[d]&
*+{\begin{ytableau}
	1&1&2\\
	2&3
\end{ytableau}}\ar@{.}_{}[d]\ar@{.}_{}[rd]&
*+{\begin{ytableau}
	1&1&1\\
	3&3
\end{ytableau}}\ar@{->}_{1}[d]\\
*+{\begin{ytableau}
	1&1&3\\
	2&3
\end{ytableau}}\ar@{->}_{2}[d]
\ar@{.}_{}[rd]&
*+{\begin{ytableau}
	1&2&2\\
	2&3
\end{ytableau}}\ar@{.}_{}[d]
&*+{\begin{ytableau}
	1&1&2\\
	3&3
\end{ytableau}}\ar@{->}_{1}[d]\\
*+{\begin{ytableau}
	1&1&3\\
	3&3
\end{ytableau}}\ar@{->}_{1}[rd]
&*+{\begin{ytableau}
	1&2&3\\
	2&3
\end{ytableau}}\ar@{.}_{}[d]&
*+{\begin{ytableau}
	1&2&2\\
	3&3
\end{ytableau}}\ar@{->}_{1}[d]\\
*{}&*+{\begin{ytableau}
	1&2&3\\
	3&3
\end{ytableau}}\ar@{->}_{1}[d]&
*+{\begin{ytableau}
	2&2&2\\
	3&3
\end{ytableau}}\ar@{.}_{}[ld]\\
*{}&*+{\begin{ytableau}
	2&2&3\\
	3&3
\end{ytableau}}
&\\
}
\end{minipage}
\end{center}
\end{center}
\vf
\end{figure}


\subsection{Kostka and 
quasi-Kostka numbers}

For a tableau $T$ of entries at most $n$, 
let 
$\wt T=(a_{1}, \ds, a_{n})$ 
where $a_{i}$ is the number of $i$ in $T$.

\begin{defn}[Kostka numbers]
For $\l$ and a weak composition $a$ 
of size $|\l|$, let 
\[
\ssyt(\l, a)=\{T\in \ssyt(\l)\mid \wt T=a\}
\qte{and}\q
K_{\l a}=|\ssyt(\l, a)|.
\]
\end{defn}

\begin{defn}[quasi-Kostka numbers]
For $\l$ and a strong composition $\a$ 
of size $|\l|$, let
\[
\st(\l, \a)=\{T\in \st(\l)\mid \des T=\al\}
\qte{and}\q
f_{\l\a}=|\st(\l, \a)|.
\]
\end{defn}
\begin{rmk}\hf
\begin{enumerate}
\item Recall the traditional notation 
$f_{\l}=|\st(\l)|$ and so $f_{\l}=\sum_{\a}f_{\l\a}$.
\item The definition of $f_{\l\a}$ implies 
\[
\sum_{T\in \F(\l, \a)}x^{\wt T}=
f_{\l\a}F_{\a}(X).
\]
\item 
Egge--Loehr--Warrington \cite{elw} and 
Orellana--Saliola--Schilling--Zabrocki \cite{ossz} 
discussed such numbers as 
entries of certain transition matrices.
\end{enumerate}
\end{rmk}

\begin{defn}
Say $T\in \ssytl$ 
is \kyo{quasi-Yamanouchi} if 
$\des T=\wt T$. 
Denote the set of all such tableaux 
by $\qy(\l)$.
\end{defn}

There is the canonical identification
$\qy(\l)\cong \st(\l)$:
standardization 
$T\mto T^{st}$, 
and destandardization 
$T\mto T^{dst}$ 
are des-preserving bijections 
and the inverse map of each other. Hence 
\[
f_{\l\a}=|\{T\in \qyl\mid \des T=\al\}|
=|\{T\in \qyl\mid \wt T=\al\}|\le 
K_{\l\a}.
\]
Observe that 
$0\le f_{\l\l}\le K_{\l\l}=1$ 
and the \kyo{superstandard tableau} $T(\l)$ 
(row $i$ consists of only $i$) 
satisfies $\des T(\l)=\wt T(\l)=\l$.
We thus showed that $f_{\l\l}=1$ for all $\l$.

\begin{ex}
SYTs of shape $\l=332$ of descent $\a=12221$ are only these:
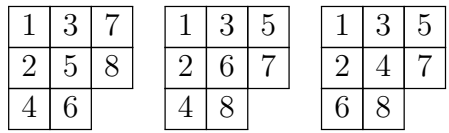
\[
\begin{ytableau}
	1&3&7\\2&5&8\\4&6
\end{ytableau}
\q
\begin{ytableau}
	1&3&5\\2&6&7\\4&8
\end{ytableau}
\q 
\begin{ytableau}
	1&3&5\\2&4&7\\6&8
\end{ytableau}
\]
Thus, $f_{\l\a}=3$.
\end{ex}


\subsection{Skeleton polynomials}

For a composition $\a$, 
denote its length by $\el(\al)$. 
In particular, for a partition $\l$, 
$\el(\l)$ is the number of rows of $\l$. 
Let $m(\l)=\max\{\el(\dest)\mid T\in \stl\}.$
\begin{lem}[{\cite[Lemma 3.1]{wa}}]
$m(\l)=n-(\l_{1}-1)$ for $\l\ne\ku$.
\end{lem}

%
By $x^{\a}$ we mean a
 monomial 
 $x_{1}^{\a_{1}}\cd x_{\el(\a)}^{\a_{\el(\a)}}$.
 
\begin{defn}
The \kyo{skeleton polynomial} for $\l$ is 
\[
\Skl(x_{1}, \ds, x_{m(\l)})=
\sum_{T\in \csl}x^{\dest}.
\]
\end{defn}
This is a homogeneous polynomial of degree 
$|\l|$ with nonnegative integer coefficients.
It has $m(\l)$ variables. 
However, if no confusion arises, we often 
abbreviate them just as $x$. 
In particular, $\Sk_{\ku}=1$.

\begin{rmk}
We may regard 
$\Sk_{\l}(x)$ as the generating function of 
$\qy(\l)$ $\k{\sub \ssyt(\l)}$:
\[
\Skl(x)=\sum_{\a\compn} f_{\l\a}x^{\a}.
\q \k
{\te{cf.}\q  
s_{\l}(X)=\sum_{\a\compn} f_{\l\a}F_{\a}(X)}\]
For this reason, we call $(f_{\l\a})$ \kyo{skeleton coefficients}.
In particular, $\Skl(1, \ds, 1)=f_{\l}$.
\end{rmk}
\begin{ex}
Table \ref{t1} shows 
$\Sk_{\l}(x)$ for $|\l|\le 4$.
\end{ex}

The cases for one-row and one-column shapes 
are 
\[
\Sk_{(c)}(x)=
x_{1}^{c}, \q 
\Sk_{(1^{r})}(x)=
x_{1}\cd x_{r}.
\]

{\renewcommand{\arraystretch}{2.25}
\begin{table}[h!]
\caption{Skeleton polynomials}
\label{t1}
\ytableausetup{mathmode, boxsize=14pt}
\begin{center}
	\begin{tabular}{cccccccc}\h
	$\l$&$\qy(\l)$	&$\Skl(x)$	\\\h
	$\ku$&$\ku$	&1	\\\h
	$\yd{1}$&	\begin{ytableau}
		1
	\end{ytableau}	&$x_{1}$
\\\h
	$\yd{2}$&	\begin{ytableau}
		1&1
	\end{ytableau}&$x_{1}^{2}$	
\\\h
	$\yd{1,1}$&	\begin{ytableau}
		1\\2
	\end{ytableau}&$x_{1}x_{2}$\vspace{3mm}
\\\h
	$\yd{3}$&	\begin{ytableau}
		1&1&1
	\end{ytableau}&$x_{1}^{3}$	\vspace{3mm}
\\\h
	$\yd{2,1}$&	\begin{ytableau}
		1&1\\2
	\end{ytableau}\q
		\begin{ytableau}
		1&2\\2
	\end{ytableau}	&$x_{1}^{2}x_{1}+x_{1}x_{2}^{2}$\vspace{3mm}
\\\h
	$\yd{1,1,1}$&	\begin{ytableau}
		1\\2\\3
	\end{ytableau}&$x_{1}x_{2}x_{3}$\vspace{3mm}\rule{0mm}{12mm}
\\\h
	$\yd{4}$&	\begin{ytableau}
		1&1&1&1
	\end{ytableau}&$x_{1}^{4}$	\vspace{2mm}
\\\h
	$\yd{3,1}$&	\begin{ytableau}
		1&1&1\\2
	\end{ytableau}\q
	\begin{ytableau}
		1&1&2\\2
	\end{ytableau}\q
	\begin{ytableau}
		1&2&2\\2
	\end{ytableau}
&$x_{1}^{3}x_{2}+x_{1}^{2}x_{2}^{2}
+x_{1}x_{2}^{3}$\vspace{3mm}	
	\\\h
	$\yd{2,2}$&	\begin{ytableau}
		1&1\\2&2
	\end{ytableau}\q
	\begin{ytableau}
		1&2\\2&3
	\end{ytableau}	&$x_{1}^{2}x_{2}^{2}
+x_{1}x_{2}^{2}x_{3}$\vspace{3mm}
	\\\h
	$\yd{2,1,1}$&	\begin{ytableau}
		1&1\\2\\3
	\end{ytableau}\q
	\begin{ytableau}
		1&2\\2\\3
	\end{ytableau}\q
	\begin{ytableau}
		1&3\\2\\3
	\end{ytableau}&$x_{1}^{2}x_{2}x_{3}+
x_{1}x_{2}^{2}x_{3}+x_{1}x_{2}x_{3}^{2}$\vspace{3mm}\rule{0mm}{12mm}
	\\\h
	$\yd{4}$&	\begin{ytableau}
		1\\2\\3\\4
	\end{ytableau}&$x_{1}x_{2}x_{3}x_{4}$	
	\\\h	\end{tabular}
\end{center}
\end{table}}


Note that $s_{\l}(X)$ is symmetric while 
$\Skl(x)$ is not necessarily.
However, as we see, it has subtle symmetry on polynomials with \kyo{distinct} number of variables.


For $\el(\l)\le i\le m(\l)$, 
define \kyo{$i$-skeleton polynomial} 
for $\l$
\[
\Sk_{\l, i}(x)=
\sum_{T\in \st(\l, \a), \el(\a)=i}
x^{\a}
\]
as an $i$-variable one. 
For all other $i$, set 
$\Sk_{\l, i}(x)=0$ for convenience.
\begin{prop}
If $\el(\l)\le i\le m(\l)$, 
then $\Sk_{\l, i}(x)\ne 0$.
\end{prop}
\begin{proof}
A crystal skeleton graph $\csl$ 
in \cite{bcds} is connected and 
if $T, U$ are adjacent vertices, 
then $\el(\des T)-\el(\des U)\in\{-1, 0, 1\}$ 
\cite[Theorem 4.26]{bcds}.
\end{proof}


%
\begin{ex}\label{ex1}
$\qy(3, 2)$ consists of 
the following five tableaux 
\[
\ytableausetup{mathmode, boxsize=15pt}
\begin{ytableau}
	1&1&1\\2&2
\end{ytableau}\q
\begin{ytableau}
	1&1&2\\2&2
\end{ytableau}
\q
\begin{ytableau}
	1&1&2\\2&3
\end{ytableau}\q
\begin{ytableau}
	1&2&2\\2&3
\end{ytableau}\q
\begin{ytableau}
	1&2&3\\2&3
\end{ytableau}
\]
so that 
$\Sk_{32, 2}(x)=x^{32}+x^{23}$ and $ 
\Sk_{32, 3}(x)=x^{221}+x^{131}+x^{122}.
$
\end{ex}

\begin{defn}
For an $i$-variable polynomial
$g(x_{1}, \ds, x_{i})$, 
define its reversal as 
\[
g^{*}(x_{1}, x_{2}, \ds, x_{i})=g(x_{i}, x_{i-1}, \ds, x_{1}).
\]
\end{defn}

\begin{thm}\label{th1}
If $\el(\l)\le i\le m(\l)$, then 
$\Sk_{\l, i}^{*}(x)=
\Sk_{\l, i}(x).$
\end{thm}
We will prove this with Theorem \ref{th2} together.

\subsection{Evacuation}

\begin{defn}[{\cite[Remark 2.5]{bcds}}]
For $T\in \B(\l)$, 
let $w=\rowt=w_{1}\cd w_{n}$ $(n=|\l|)$ be its row word (reading from bottom to top, left to right).
Let $w^{*}=(n+1-w_{n})\cd (n+1-w_{1})$ 
and define the \kyo{evacuation} (Lusztig/Sch\"{u}tzenberger involution) of $T$ as 
$T^{*}=P(w^{*})$. 
\end{defn}

\begin{ex}
Observe that 
\[
T=
\begin{ytableau}
	1&1&1&2\\
	3&4\\4
\end{ytableau}, 
w=\row T=4341112, 
w^{*}=6777454, 
T^{*}=
P(w^{*})=
\begin{ytableau}
	4&4&7&7\\
	5&7\\6
\end{ytableau}
\]
with $\des T=412$, 
$\des (T^{*})=214$. 
\end{ex}

For 
$\a=(\a_{1}, \ds, \a_{l})\compn$, let 
$\a^{*}=(\a_{l}, \ds, \a_{1})$ denote its reversal.

\begin{lem}[See {\cite[Section 4.4]{bcds}}]
Let $T \in \B(\l)$. 
\begin{enumerate}
\item $\des (T^{*})=(\dest)^{*}$. 
	In particular, 
	$\el(\des (T^{*}))=\el(\dest)$.
\item 
Let $\B_{N}(\l)$ (for $N\ge |\l|$) be a bounded crystal graph with all SSYTs of entries at most $N$. Then $*:\B_{n}(\l)\to \B_{n}(\l)$ $(n=|\l|)$
is an involutive crystal anti-isomorphism, that is, 
\[
({T^{*}})^{*}=T
\qte{and}\q
T\xto[i]{}U
\iff
U^{*}\xto[n-i]{}T^{*} \q(1\le i\le n-1)
\] 
where $T\xto[i]{}U$ means $U=\f_{i}(T)$ 
and $T, U\ne 0$. 
In particular, it restricts to an involutive bijection on $\stl$.
\end{enumerate}
\end{lem}


\begin{thm}
\label{th2}
For all $\l$, $\a$, we have 
$f_{\l\a}=f_{\l\a^{*}}$.
As a consequence, Theorem \ref{th1} holds.
\end{thm}

\subsection{Inner crystal}

Note that $\el(\l)=\min\{\el(\dest)\mid T\in \stl\}$ 
because the superstandard tableau $T(\l)$ satisfies $\el(\des T(\l))=\el(\l)$ and clearly each 
$T\in \B(\l)$ contains at least 
$\el(\l)$ maximal horizontal bands.

\begin{defn}
Define the \kyo{inner crystal} of $\B(\l)$ as 
\[
\A(\l)=
\{T\in \csl\mid \el(\des T)=\el(\l)\}.
\]
\end{defn}
There is a specific reason to call this a "crystal".



\begin{thm}
After arranging vertices and edges appropriately, $\A(\l)$ 
is isomorphic to $\B_{\el(\l)}(\l)$ as a crystal graph.
\end{thm}
\begin{proof}
See \cite[Theorem 4.37]{bcds} for details.
\end{proof}

\begin{rmk}
We can say that the inner crystal of $\B_{N}(\l)$ 
is also $\B_{\el(\l)}(\l)$. Inclusions 
\[
\A(\l)=\B_{\el(\l)}(\l)\sub \csl \sub \B_{N}(\l)\sub \B(\l)
\]
indicate ``saturated sub-crystal property" of  ambient crystals. 

%
%
%
\end{rmk}


\begin{defn}
The \kyo{inner crystal polynomial} for $\l$ is 
$A_{\l}(x)=\Sk_{\l, \el(\l)}(x)$.
\end{defn}
\begin{thm}
$A_{\l}(x)=s_{\l}(x_{1}, \ds, x_{\el(\l)}).$
In particular, this is symmetric.
\end{thm}
\begin{ex}
\ytableausetup{mathmode, boxsize=8pt}
$\B_{3}\k{\yd{2,1}}$
consists of two quasi-crystals.
\[
\ytableausetup{mathmode, boxsize=15pt}
\begin{ytableau}
	1&1\\2
\end{ytableau}
\too
\begin{ytableau}
	1&1\\3
\end{ytableau}
\too
{\bf{\begin{ytableau}
	\bf{1}&\bf{2}\\\bf{3}
\end{ytableau}}}
\too
\begin{ytableau}
	2&2\\3
\end{ytableau}
\]
\[
\begin{ytableau}
	1&2\\2
\end{ytableau}
\too
{\bf{\begin{ytableau}
	\bf{1}&\bf{3}\\\bf{2}
\end{ytableau}}}
\too
\begin{ytableau}
	1&3\\3
\end{ytableau}
\too
\begin{ytableau}
	2&3\\3
\end{ytableau}
\]
Two standard tableaux forms its skeleton and the inner crystal. It is easier to see the edge between them after destandardization:
\[
{\bf{\begin{ytableau}
	1&1\\2
\end{ytableau}}}
\q 
\underset{1}{\too}
\q
{\bf{\begin{ytableau}
	1&2\\2
\end{ytableau}}}
\q 
\ytableausetup{mathmode, boxsize=8pt}
\k{\cong \B_{2}\k{\yd{2,1}}}
\]
\end{ex}


\section{Calculus of Skeleton polynomials}

Rather than studying a solo skeleton polynomial, 
it is desirable to investigate interactions 
over a family of those. We confirm several basic facts in this section.

\subsection{Linear indepenency}

For $\l, \m\in \te{Par}(n)$, 
define the \kyo{dominance order} 
$\l\ge \mu$ if 
\[
\l_{1}+\cd+\l_{i}\ge \m_{1}+\cd+\m_{i}
\]
for all $i$ (adding 0's to the end if necessary). Extend this order onto 
$\te{Comp}(n)$, the set of all 
strong compositions of size $n$.

\begin{lem}\label{l2}
If $\l\not\ge \a$, then $f_{\l\a}=0$.
As a consequence, 
\[
\Skl(x)=x^{\l}+\sum_{\a<\l}f_{\l\a}x^{\a}.
\]
\end{lem}

\begin{proof}
This follows from 
$f_{\l\l}=1$, $0\le f_{\l\a}\le K_{\l\a}$ 
and $K_{\l\a}=0$ if $\l\not\ge \a$.
\end{proof}



\begin{thm}
$(\Sk_{\l}(x)\mid \l\parn)$ is linearly independent 
in $\cc[x_{1}, \ds, x_{n}]$.
In particular, 
$\Sk_{\l}(x)=
\Sk_{\m}(x)\then \l=\m$.
\end{thm}
\begin{proof}
Let $N=|\Parn|$
and choose a sequence 
$\l^{(1)}, \ds, 
\l^{(N)}$ from $\Parn$ such that 
$\l^{(i)}$ is maximal 
in the subposet 
$\{\l^{(i)}, \ds, \l^{{(N)}}\}$.
To prove the 
linear independency, 
suppose 
\[
\dsum_{i=1}^{N}c_{i}\Sk_{{\l^{(i)}}}(x)=0, \q c_{i}\in \cc.
\]
Thanks to Lemma \ref{l2}, the coefficient of 
$x^{\l^{(1)}}$ comes from 
only $\Sk_{{\l^{(1)}}}(x)$, 
that is, 
\[
c_{1}
x^{\l^{(1)}}+
\k
{\te{
all other terms without $x^{\l^{(1)}}$}}=0.
\]
Thus, $c_{1}$ must be $0$.
Similarly, we can show $c_{2}=\cd=c_{N}=0$.
\end{proof}


\subsection{Skeleton RSK correspondence}

Recall the classical \kyo{RSK correspondence} 
\[
\K_{n}\to \bigsqcup_{\l\parn}\ssyt(\l)\ti \ssyt(\l), 
\q w\mto(P(w), Q(w))
\]
where $\K_{n}$ is 
the set of all Knuth arrays (lexicographic biwords) of length $n$. 
According to Stanley \cite[p.419]{sta}, we have 
\[
\sum_{\l\,\parn}s_{\l}(X)s_{\l}(Y)=
\sum_{w\in S_{n}} 
F_{\des w^{-1}}(X)F_{\des w}(Y)
\]
as $\des\, w^{-1}=\des P(w), 
\des\, w=\des Q(w)$. 
Moreover, it restricts to 
each of the following bijections
\begin{align*}
	S_{n}^{\te{invol}}&\to \bigsqcup_{\l\parn}\stl,\q w\mto P(w),\\
	S_{n}&\to \bigsqcup_{\l\parn}\stl\ti \stl, \q w\mto(P(w), Q(w)),
	\\
	\nn^{n}&\to \bigsqcup_{\l\parn}\ssyt(\l)\ti \stl, \q w\mto(P(w), Q(w))
\end{align*}
where 
$S_{n}^{\te{invol}}$
 is the set of all involutions in $S_{n}$, 
$[k]^{n}$ the set of all words 
of $\{1, \ds, k\}$ of length $n$ 
and 
$\nn^{n}$ the set of all words of 
length $n$. 
 

\begin{thm}[Skeleton R, RS correspondences]\label{tsk}
\begin{align*}
	\sum_{\l\partition n}\Sk_{\l}(x)&=\sum_{w\in S_{n}^{\te{invol}}}
x^{\des w}.
	\\\sum_{\l\partition n}\Sk_{\l}(x)\Sk_{\l}(y)
	&=\sum_{w\in S_{n}}x^{\des(w^{-1})}
y^{\des w}.
	\\\sum_{\l\partition n}
s_{\l}(X)\Sk_{\l}(y)
&=\sum_{w\in S_{n}}F_{\des(w^{-1})}
(X)y^{\des w}.
\end{align*}
\end{thm}

\begin{proof}
Each proof of these is quite similar.
Here we prove only the second identity.
\[
\te{RHS}=
\sum_{\l\parn}
\sum_{(S, T)\in \st(\l)\ti \st(\l)}
x^{\des S}y^{\des T}
=
\sum_{\l\parn}
\sum_{(S^{\dst}, T^{\dst})\in \qy(\l)\ti \qy(\l)}
x^{\des S^{dst}}y^{\des T^{dst}}
\]
\[
=
\sum_{\l\parn}
\k
{\sum_{U\in \qy(\l)}x^{\des U}}
\k{\sum_{V\in \qy(\l)}y^{\des V}}
=\te{LHS}.
\]
\end{proof}

\begin{ex}
Table \ref{t1} shows that 
\[
\sum_{\l\partition 3}\Sk_{\l}(x)\Sk_{\l}(y)=
x_{1}^{3}y_{1}^{3}+
(x_{1}^{2}x_{2}+x_{1}x_{2}^{2})
(y_{1}^{2}y_{2}+y_{1}y_{2}^{2})
+x_{1}x_{2}x_{3}y_{1}y_{2}y_{3}.
\]
\end{ex}


\subsection{Skeleton enumeration}

Further, we can use skeleton polynomials to 
enumerate permutations with restrictions 
on length of descent.

\begin{cor}
Let $\el(\l)\le i\le m(\l)$. Then 
\[
\sum_{\l\parn}\Sk_{\l}(1^{i}, 0^{m(\l)-i})
\]
equals the number of involutions $w\in S_{n}$ such that 
$\el(\des (w))\le i$.
In particular, 
\[
\Sk_{\l}(1^{m(\l)})=
\sum_{\l\,\partition n}f_{\l}
=|S_{n}^{\te{invol}}|.
\]
\end{cor}
\begin{ex}
Table \ref{t1} shows that 
\[
\Sk_{4}(1)
+\Sk_{31}(1, 1)+\Sk_{22}(1, 1, 0)+
\Sk_{211}(1, 1, 0)+
\Sk_{1111}(1, 1, 0, 0)
\]
\[=1+3+(1+0)+0+0=5
\]
counts five involutions $w\in S_{4}$
 with $\el(\des w)\le 2$ as 
1234, 1243, 1324, 2134, 3412.
\end{ex}


\begin{cor}
Let $\el(\l)\le i, j \le m(\l)$. Then 
\[
\sum_{\l\parn}\Sk_{\l}(1^{i}, 0^{m(\l)-i})
\Sk_{\l}(1^{j}, 0^{m(\l)-j})
\]
equals the number of permutations 
$w\in S_{n}$ such that 
$\el(\des(w^{-1}))
\le i$ and $\el(\des w)\le j.$
In particular, 
$\sum_{\l}f_{\l}^{2}=n!$.
\end{cor}

\subsection{Skeleton hook sum}

A nonempty partition $\l$ is a \kyo{hook} if 
$\l_{2}\le 1$ (i.e., $\l=(\l_{1}, 1^{\el(\l)-1})$).
This definition includes one-row and one-column shapes.
\begin{lem}
$\l$ is a hook $\iff \el(\l)=m(\l)$
$\iff$ $\Sk_{\l}(x)=\Sk_{\l, \el(\l)}(x)$.
\end{lem}

\begin{thm}
$\sum_{\al\compn}x^{\a}=
\sum_{
\l\parn,\, \te{$\l$ a hook}
}\Sk_{\l}(x).$
\end{thm}
\begin{proof}
For each $\al\compn$ of length $k$ $(1\le k\le n)$, there exists a unique hook of size $n$ of length $k$, that is, $\l^{(k)}=(n-k+1, 1^{k-1})$, 
such that $|\st(\l^{(k)}, \a)|=1$. 
\end{proof}

\begin{ex}
\ytableausetup{mathmode, boxsize=5pt}
\[
\sum_{\al\models 4}x^{\a}=
\Sk_{\,\yd{4}}(x_{1})+
\Sk_{\,\yd{3,1}}(x_{1}, x_{2})+
\Sk_{\,\yd{2,1,1}}(x_{1}, x_{2}, x_{3})+
\Sk_{\,\yd{1,1,1,1}}(x_{1}, x_{2}, x_{3}, x_{4}).
\]
\end{ex}

\begin{figure}
\caption{$(\te{Comp}(4), \le, \dep)$}
\label{figcomp}
\begin{center}
\mb{}\\
\begin{minipage}[c]{.8\tw}
\xymatrix@=5mm
{
*{\te{depth}}&0&1&2&3&4&5&6\\
&\mb{}&&&*+{211}\ar@{-}[rd]&&&\\
&*+{4}\ar@{-}[r]&*+{31}\ar@{-}[r]&*+{22}\ar@{-}[ru]\ar@{-}[rd]&&*+{121}\ar@{-}[r]&*+{112}\ar@{-}[r]&*+{1111}\\
&\mb{}&&&*+{13}\ar@{-}[ru]&&&\\
}
\end{minipage}
\end{center}
\end{figure}

\begin{figure}
\caption{$(2^{[3]}, \le, \maj)$}
\label{fig1}
\begin{center}
\mb{}\\
\begin{minipage}[c]{.8\tw}
\xymatrix@=5mm
{
*{\te{maj}}&0&1&2&3&4&5&6\\
&\mb{}&&&*+{\{1, 2\}}\ar@{-}[rd]&&&\\
&*+{\ku}\ar@{-}[r]&*+{\{1\}}\ar@{-}[r]&*+{\{2\}}\ar@{-}[ru]\ar@{-}[rd]&&*+{\{1, 3\}}\ar@{-}[r]&*+{\{2, 3\}}\ar@{-}[r]&*+{\{1, 2, 3\}}\\
&\mb{}&&&*+{\{3\}}\ar@{-}[ru]&&&\\
}
\end{minipage}
\end{center}
\end{figure}

\section{Deep calculus}

\subsection{Depth}

Finally, we develop ``deep calculus" of 
skeleton polynomials and Schur functions with 
the following statistic. 
%
For $\a\models n$, its \kyo{depth} is 
$\ch(\a)=\sum_{i=1}^{\el(\a)}
(i-1)\a_{i}.$ This is meaningful in our context 
with $(f_{\l\a}|\l\ge \a)$ as follows.
\begin{prop}\label{p1}
The poset $(\te{Comp}(n), \le, \dep)$ 
is graded with lowest depth 0 and 
highest $\binom{n}{2}$.
\end{prop}
\begin{proof}
Clearly, $\dep(n)=0$ and
$\dep(1^{n})=\binom{n}{2}$ and only such 
compositions take the extreme depth.
Any $\al\in \te{Comp}(n)$ can be obtained from 
$(n)$ by a sequence of 
the following raising/refinement operations:
\begin{align*}
	R_{i}\a&=(\a_{1}, \ds, \a_{i-1}, \a_{i}-1, \a_{i+1}+1, \a_{i+2}, \ds, \al_{\el(\a)})
\qte{if\q} \a_{i}\ge 2, \,\,i\le \el(\a)-1,
	\\R_{\el(\a)}\a&=(\a_{1}, \ds, \al_{\el(\a)-1}, \al_{\el(\a)}-1, 1)
\qte{if\q} \a_{\el(\a)}\ge 2.
\end{align*}
Check that if $\al\ne (1^{n})$, 
there exists at least one $i$ such that 
$R_{i}\a$ is defined, $\al \rhd R_{i}\a$ (a covering relation in dominance order) and 
$\dep(R_{i}\a)-\dep(\a)=1$.
\end{proof}
Figure \ref{figcomp} shows an example of such posets. 
Now for $T\in \ssytl$, define its \kyo{depth} as 
$\ch(T)=\ch(\dest)$. For example,
\ytableausetup{mathmode, boxsize=14pt} 
\[
\dep\k{\,
\begin{ytableau}
	1&1&4&4\\
	2&4\\3
\end{ytableau}\,}=\dep(214)=9.
\]
It is now natural to think of 
\kyo{deep Schur function}
\[
s_{\l}(X, q)=\sum_{T\in \ssytl}q^{\ch (T)}x^{\wt (T)}.
\]
Indeed, $\ch(T)$ depends on only $\des T$ 
so that, as Assaf mentioned \cite{as}, it really makes sense to talk about 
deep Gessel expansion and 
skeleton polynomials
\[
s_{\l}(X, q)=
\sum_{\a}f_{\l\a}q^{\ch(\a)}F_{\a}(X), 
\q
\Sk_{\l}(x, q)=
\sum_{\a}f_{\l\a}q^{\ch(\a)}x^{\a}.
\]
\ytableausetup{mathmode, boxsize=5pt}
Observe that 
Example \ref{ex1} extends to 
\[
\Sk_{\,\yd{3,2}}(x, q)=
q^{2}x^{32}+q^{3}x^{23}+q^{4}x^{221}+q^{5}x^{131}+q^{6}x^{122}.
\]

%

\subsection{Major index}

Now is time to explain technical twists among 
depth, \kyo{descent sets} and \kyo{major index}.


Consider two sets $\te{Comp}(n)$ and $2^{[n-1]}$ of both cardinality $2^{n-1}$. 
There is the traditional identification:
For 
$\a=(\a_{1}, \ds, \a_{l})\in \te{Comp}(n)$, 
define 
$A(\a)=\{a_{1}, \ds, a_{l-1}\}$ with $a_{i}=\a_{1}+\cd+\a_{i}$. On the other hand, 
given 
$A=\{a_{1}, \ds, a_{k}\}\in 2^{[n-1]}$, say 
$a_{1}<\cd <a_{k}$ $(0\le k\le n-1$), 
define 
$\a(A)=(\a_{1}, \ds, \a_{k+1})$ 
with 
$\a_{1}=a_{1}, 
\a_{i}=a_{i}-a_{i-1}$
$(2\le i\le k)$ and 
$\a_{k+1}=n-a_{k}$.
In fact, 
$\a\mto A(\a), 
A\mto \a(A)$
are inverse maps of each other. 
To go into more details, introduce a superboolean structure onto $2^{[n-1]}$: for $A, B\in 2^{[n-1]}$, 
declare $A\lhd B$ if either 
\begin{itemize}
	\item $A\subset B, 1\not\in A, 1\in B$, 
	and $B\sm \{1\}=A$ or 
	\item 
	there exists some $a\in A$ such that 
	$a\not\in B$,
	$a+1\not\in A$,$a+1\in B$
	and $B\sm\{a+1\}=A\sm\{a\}$.
\end{itemize}
Define $A\le B$ by the transitive closure of this relation so that $2^{[n-1]}$ forms a poset graded by $\maj(A)=\sum_{a\in A}a$ by construction. 
It is isomorphic to $(\te{Comp}(n), \le, \dep)$ as sets.
However, maps $\a\mto A(\a), 
A\mto \a(A)$ are \kyo{not} order-preserving nor order-reversing.
For example, $\dep(214)=9$ $(n=7, \binom{7}{2}=21)$ while $\maj\{2, 3\}=5$.
Let us clarify the precise correspondence between depth and major index.
For $T\in \stl$, 
let 
\[
\Des\, T=\{i\in [n-1]\mid \te{$i+1$ is strictly lower than $i$ in $T$}\}.
\]
More generally, for $T\in \ssytl$, 
define $\Des\, T=\Des (T^{st})$ and  
$\maj(T)=\sum_{i\in \Des T}i$.
Observe that 
$\des T=(\a_{1}, \ds, \a_{l})
\iff
\Des \,T=\{\a_{1}, \a_{1}+\a_{2}, \ds, 
\a_{1}+\cd+\a_{l-1}\}.$
Note that $\des T$ and $\maj\,T$ are \kyo{not}  termwise equal.

\begin{lem}\label{lmaj}
For each $T\in \ssytl$, 
$\dep(T)=\maj(T^{*})$. 
\end{lem}
\begin{proof}
Say $\des T=\al=(\a_{1}, \ds, \a_{l})\compn$.
Then 
$\des (T^{*})=
(\des T)^{*}=
(\a_{l}, \ds, \a_{1}).$
\[
\Des (T^{*})=
\{\a_{l}, \a_{l}+\a_{l-1}, \ds, \al_{l}+\cd+\a_{2}\}.
\]
\[
\maj \k{\Des (T^{*})}=
\a_{l}+
(\a_{l}+\a_{l-1})+(\al_{l}+\cd+\a_{2})
=\sum (i-1)\a_{i}=\dep\,\al.
\]
This means $\dep(T)=\maj(T^{*})$. 
\end{proof}

There is another application of 
Lemma \ref{lmaj} to symmetric groups 
and $q$- or $(p, q)$-analogs of 
all identities in Theorem \ref{tsk}. 

For $w\in S_{n}$, 
define $\des w=\des Q(w)$, $\dep\, w=\dep\, Q(w)$.

\begin{thm}
\label{tskq}
\begin{align*}
\sum_{\l\partition n}\Sk_{\l}(x, p)
&=\sum_{w\in S_{n}^{\te{invol}}}
x^{\des w}p^{\dep\, w}.
\\
\sum_{\l\partition n}\Sk_{\l}(x, p)\Sk_{\l}(y, q)
&=\sum_{w\in S_{n}}x^{\des(w^{-1})}
p^{\dep\, w^{-1}}
y^{\des w}q^{\dep\, w}.
\\
\sum_{\l\partition n}
s_{\l}(X)\Sk_{\l}(y, q)
&=\sum_{w\in S_{n}}F_{\des(w^{-1})}
(X)y^{\des w}q^{\dep\, w}.
\end{align*}
\end{thm}
\begin{rmk}
$\sum_{\l\vdash n}s_{\l}(X)\Sk_{\l}(y, q)$ 
with $y_{i}=1$ recovers 
the graded Frobenius character of $S_{n}$.
\end{rmk}

Define 
$\maj\, w=\maj\, Q(w)$. 
As is well-known, this is a Mahonian statistic.

\begin{thm}
Depth is Mahonian.
\end{thm}
\begin{proof}
Recall that 
$w^{-1}\leftrightarrow w \leftrightarrow w^{*}$ are bijections on $S_{n}$ and 
$P(w^{-1})=Q(w)$ by Symmetry Theorem on RSK.
\[
\sum_{w\in S_{n}}q^{\depw}
=\sum_{}q^{\dep\, Q(w)}
=\sum_{}q^{\dep\, P(w^{-1})}
=\sum_{}q^{\dep\, P(w^{*})}
\]
\[=\sum_{}q^{\dep\, P(w)^{*}}
=\sum_{}q^{\maj\, P(w)}
=\sum_{}q^{\maj\, Q(w^{-1})}
=\sum_{}q^{\maj\, Q(w)}
=\sum_{}q^{\maj\, w}=[n]_{q}!.
\]
\end{proof}


\subsection{Fake degree polynomials}

%

\begin{defn}[Stanley \cite{sta79}]
The \kyo{fake degree polynomial} for $\l$ is 
\[
f_{\l}(q)=\sum_{T\in\stl}q^{\maj (T)}.
\]
\end{defn}

These polynomials play an important role in representation theory and 
invariant theory. Billey--Konvalinka--Swanson \cite{bks} studied them from the following viewpoint. Say a polynomial $g(q)=\sum_{i=0}^{N}a_{i}q^{i}$ has an \kyo{internal zero} if $a_{j}\ne 0, a_{k}=0$ and $a_{l}\ne 0$ for some $j<k<l$. 
Here we count it in this way:
\[
\#\{k\in [N]\mid \te{$a_{j}\ne 0, a_{k}=0$ and $a_{l}\ne 0$ for some $0\le j<k<l\le N$}\}
\]
They proved that 
$f_{\l}(q)$ has so rarely (but possibly) an internal zero. Before stating their theorem, 
let us introduce these words for convenience.
\begin{defn}
Say $\l$ is \kyo{regular} if 
$\l\ne (c^{r})$ for any integers $c, r\ge 2$. 
Otherwise, it is \kyo{irregular}.
\end{defn}

\begin{thm}[{\cite[Theorem 1.1]{bks}}]
\label{tbks}\hf
\begin{enumerate}
	\item If $\l$ is regular, then 
$f_{\l}(q)$ has no internal zero.
Moreover, 
\[
[q^{\dep(\l)}](f_{\l})=
[q^{\binom{n}{2}-\dep(\l')}](f_{\l})=1
\]
and these are coefficients for terms 
of lowest and highest degree.
	\item Suppose $\l$ is irregular and 
$0\le k\le \binom{n}{2}$.
Then 
\[
[q^{k}](f_{\l})=0
\]
if and only if 
$k<\dep (\l)$ or 
$k=\dep(\l)+1$ 
or $k=\binom{n}{2}-\dep(\l')-1$ 
or $k>\binom{n}{2}-\dep(\l')$.
As a consequence, 
$f_{\l}(q)$ has at least one and at most two  internal zero(s).
Moreover, 
\[
[q^{\dep(\l)}](f_{\l})=
[q^{\binom{n}{2}-\dep(\l')}](f_{\l})=1.
\]
\end{enumerate}

\end{thm}



They devoted the full Section 4 in \cite{bks} to the proof. Here, we give a short proof for 
a part of these assertions.

\begin{lem}\label{lbks}
$f_{\l}(q)=\sum_{T\in\stl}q^{\dep (T)}.$
\end{lem}
\begin{proof}
Recall that $T\leftrightarrow T^{*}$ (evacuation) is an involutive bijection on $\stl$. 
It follows from Lemma \ref{lmaj} that 
\[
f_{\l}(q)=
\sum_{T\in \stl}q^{\maj(T)}
=\sum_{T\in \stl}q^{\dep(T^{*})}
=\sum_{T\in \stl}q^{\dep(T)}.
%
\]
\end{proof}

\begin{ex}
Observe from Table \ref{t1} that 
\ytableausetup{mathmode, boxsize=4pt}
$f_{\,\yd{4}}(q)=1$,
$f_{\,\yd{3,1}}(q)=q+q^{2}+q^{3}$,
$f_{\,\yd{2,2}}(q)=q^{2}+q^{4}$,
$f_{\,\yd{2,1,1}}(q)=q^{3}+q^{4}+q^{5}$,
$f_{\,\yd{1,1,1,1}}(q)=q^{6}$.
\end{ex}

\begin{thm}
Let $\l$ and $k$ as in Theorem \ref{tbks} (2). 
If $k<\dep (\l)$ or 
$k=\dep(\l)+1$, then 
$[q^{k}](f_{\l})=0.$
\end{thm}
\begin{proof}
Say $\l=(c^{\el(\l)}), c, \el(\l)\ge 2$.
Suppose 
$[q^{k}](f_{\l})\ne 0$, 
that is, there exists some $T\in\stl$ such that 
$\dep(T)=k$. 
We then must have $\l \ge \des T$ in dominance order.
Therefore $\dep(\l)\le \dep (\des T)=\dep\, T=k$
 forces $k\not<\dep (\l)$.
Now if $k=\dep(\l)+1$, then $\l$ covers $\des T$ in $\te{Comp}(n)$. 
According to Proposition \ref{p1}, there are two cases.
\begin{enumerate}
\item $R_{i}\l=\des T$ for some $i\le \el(\l)-1$.
Then 
$\des T=(c^{i-1}, c+1, c-1, c, \ds, c)$.
However, the size of a horizontal strip in $T$ on  shape $(c^{\el(\l)})$ must be at most $c$. 
This is impossible.
\item $R_{\el(\l)}\l=\dest$, 
i.e., 
$\dest=(c^{\el(\l)-1}, c-1, 1)$.
This means the first $\el(\l)-1$ rows 
of $T$ coincides with ones of the superstandard tableau $T(\l)$. However, 
there cannot exist two more 
maximal horizontal bands in the last one row.
\end{enumerate}
\end{proof}



\subsection{Schur family}



\begin{defn}
For $\l\parn$, define 
the \kyo{supersemistandard tableau} 
$Q(\l)\in \qyl$ with 
$\des Q(\l)=\l$. 
Let $\l'=(\l_{1}', \ds, \l_{l}', 1^{L})$ $(\l_{l}'\ge 2, L\ge0)$ denote the conjugate of $\l$.
The \kyo{anti-supersemistandard tableau} 
$\ol{Q}(\l)$ is 
\[
\begin{array}{ccccccc}
1&\l_{1}'   &\ds
&\l_{1}'+\cd+\l_{l-1}'-(l-2)&m(\l)&\cd &m(\l)   \\
2&\l_{1}'+1  &\cd
&
\l_{1}'+\cd+\l_{l-1}'-(l-1)\\
\vdots&\vdots&&\vdots\\
	\vdots&\vdots&\cd&	
	m(\l)\\
	\vdots&
	\l_{1}'+\l_{2}'-1&&\\
	\l_{1}'&\\
\end{array}
\]
Let 
$\oll=\des \ol{Q}(\l)
=(1^{\l_{1}'-1}, 2, 1^{\l_{2}'-2}, 2, \ds, 2, 1^{\l_{l}'-2}, L+1)\,(\in \te{Comp}(n))$.
\end{defn}
By construction, 
$Q(\l)$ $(\ol{Q}(\l))$ is the unique tableau taking minimal (maximal) depth in $\qyl$ as its existence guaranteed by Theorem \ref{tbks} and Lemma \ref{lbks}. 
In particular, $\dep (\oll)=
\binom{n}{2}-\dep(\l')$. 
\begin{rmk}
If $\l$ is a one-row or one-column shape, 
then $\ol{Q}(\l)=Q(\l), \ol{\l}=\l$.
\end{rmk}

\begin{ex}
For $\l=331$, we have 
$\l'=322$, 
\[
\ytableausetup{mathmode, boxsize=14pt}
{Q}(\l)=
\begin{ytableau}
	1&1&1\\
	2&2&2\\
	3
\end{ytableau}, \q
\ol{Q}(\l)=
\begin{ytableau}
	1&3&4\\
	2&4&5\\
	3
\end{ytableau}, \q
\ol{\l}=11221.
\]
\end{ex}

\renewcommand{\g}{\gamma}
\renewcommand{\d}{\delta}

\begin{lem}\hf
\begin{enumerate}
\item 
Let $n\ge 2$, $\a, \g\models n$ and $\b, \d\models (n-1)$. 
Assume that 
\[
\b=(\b_{1}, \ds, \b_{k}), 
\b_{k}\ge 0, 
\a=(\b_{1}, \ds, \b_{k-1}, \b_{k}+1),
\]
\[
\d=(\d_{1}, \ds, \d_{m}), k\le m
\]
and 
\[
\g=(\d_{1}, \ds, \d_{m-1}, \d_{m}+1)
\]
or 
\[
\g=(\d_{1}, \ds, \d_{j}, 1, \d_{j+1}, \ds, \d_{m}).
\]
Then, 
$\b\ge \d$ implies 
$\a\ge \g$.
\item If $\l\parn, \a\models n$ and $f_{\l\a}\ne0$, then 
$\a\ge \oll.$
\end{enumerate}
\end{lem}

\begin{proof}\hf
\begin{enumerate}
\item Let $\a, \b, \g, \d$ as assumed 
and suppose $\b\ge \d$. 
For convenience, set 
$a_{i}=(\a_{1}+\cd+\a_{i})-
(\g_{1}+\cd+\g_{i})$. Note that anyway
\[
a_{k}=n-
(\g_{1}+\cd+\g_{k})\ge 0.
\]	
Case 1: if 
$\g=(\d_{1}, \ds, \d_{m-1}, \d_{m}+1)$, 
then for $i<k$,
\[
a_{i}=(\b_{1}+\cd+\b_{i})-
(\d_{1}+\cd+\d_{i})\ge 0.
\]
Case 2: if 
$\g=(\d_{1}, \ds, \d_{j}, 1, \d_{j+1}, \ds, \d_{m})$
, then for $i\le j$, 
clearly $a_{i}\ge 0$
and for $i$ such that 
$j+1\le i<k$, 
\[
a_{i}=
(\a_{1}+\cd+\a_{i})-
(\g_{1}+\cd+\g_{i})
\ge 
(\b_{1}+\cd+\b_{i})-
(\d_{1}+\cd+\d_{i-1}+1)
\]
\[=(\b_{1}+\cd+\b_{i-1})
-(\d_{1}+\cd+\d_{i-1})
+(\b_{i}-1)\ge0.
\]
At any case, we showed that 
$\a\ge \g$.
\item 
Induction on $n\ge 1$. 
If $n=1$, then 
clearly $\a=(1)=\oll$. 
Suppose $n\ge 2$ and 
$f_{\l\a}\ne 0$. 
There exists some $T\in \qyl$ such that 
\[
\des T=\al, \q\te{say } \a=(\a_{1}, \ds, \a_{k}).
\]
Now let $U$ be the (necessarily quasi-Yamanouchi) 
tableau obtained from $T$ by removing the rightmost largest entry. 
Say 
$\m=\te{shape}\, U\vdash(n-1)$ and 
$\b=\des U\models(n-1)$. 
We can write 
\[
\b=(\b_{1}, \ds, \b_{k-1}, \b_{k}), \b_{k}\ge 0
\]
and so 
\[
\a=(\b_{1}, \ds, \b_{k-1}, \b_{k}+1).
\]
Comparing $\ol{Q}(\l)$ and 
$\ol{Q}(\m)$, 
$\ol{\mu}$ is obtained from $\oll$ by 
removing one 1 or replacing $L+1$, the last entry of $\oll$, by $L$.
In other words, 
$\ol{\l}$ is obtained from $\ol{\mu}$ by 
inserting one 1 or adding 1 to the last part.
Inductive hypothesis assures now that 
$\b\ge \ol{\mu}$. 
Thanks to (1) with 
$\delta=\ol{\mu}, 
\g=\ol{\l}$, we conclude that 
$\alpha\ge \g=\ol{\l}$. 
\end{enumerate}
\end{proof}

Let 
\[
\lef[\b, \g\ri]=\{\a\in \te{Comp}(n)
\mid \gam\ge \al \ge \b\}
\]
be a symbol for an interval. 
With all our discussions, 
we have proved the following theorem.

\begin{thm}\label{treg}
If $\l$ is regular, then 
\[
s_{\l}(X, q)=\sum_{\a}f_{\l\a}q^{\ch(\a)}F_{\a}(X) 
\]
 has no internal zero. 
Moreover, $\te{Sch}(\l)\sub \lef[\ol{\l}, \l\ri]$ 
and $f_{\l\l}=f_{\l\oll}=1$.
\end{thm}

Let $C\sub \te{Comp}(n)$.  
Say it is \kyo{connected} if 
the induced subgraph with vertices $C$ in Hasse diagram of $(\te{Comp}(n), \le)$ is connected as an undirected graph.




\begin{thm}\label{tirreg}
If $\l$ is irregular, then 
$s_{\l}(X, q)$ has at least one 
and at most two internal zero(s).
Moreover, 
$\te{Sch}(\l)\sub \left[\ol{\l}, \l\ri]$ 
and it is disconnected.
\end{thm}

\begin{cor}
Suppose $\l$ is irregular. Then 
$f_{\l}(q)$ has exactly one internal zero 
$\iff \l=(2, 2)$.
\end{cor}

\begin{proof}
Suppose $\l$ is irregular. 
$f_{\l}(q)$ has exactly one internal zero 
$\iff$ it has only two nonzero terms 
$\iff$ $|\stl|=2$
$\iff \l=(2, 2)$.
\end{proof}

Comparing Theorems \ref{treg} and \ref{tirreg}, 
we come to an interested question.

\begin{prob}
Suppose $\l$ is regular. 
When is $\te{Sch}(\l)$ connected?
\end{prob}

%
%



%


\subsection{Charge, depth, inversion}

In addition to major index, 
we again describe a technical twist on depth and \kyo{charge} as Lascoux-Sch\"{u}tzenberger introduced \cite{ls}.

We are going to 
define $c(w)$ for $w\in S_{n}$ with inductive 
labelings $c_{i}(w)\in \{0, \ds, i-1\}$.
For convenience, 
introduce the set of \kyo{left descents}
\[
D_{L}(w)=\{i\in [n-1]\mid 
w^{-1}(i+1)<w^{-1}(i)\}
\]
as in the context of Coxeter groups. 
First, set $c_{1}(w)=0$.
For $i\ge2$, let 
$c_{i}(w)=c_{i-1}(w)$ if 
$i-1\not\in D_{L}(w)$.
Otherwise, $c_{i}(w)=c_{i-1}(w)+1$.
Continue this up to $i=n$. 
Define $c(w)=\sum_{i=1}^{n}c_{i}(w)$ 
$\in\{0, \ds, \binom{n}{2}\}$.

\begin{rmk}
Some authors call our $c(w)$ ``cocharge". 
We have to be careful whenever we deal with such a statistic in this context.
\end{rmk}

\begin{ex}
Let $w=57841362$. 
Labelings are 
\[
\begin{array}{c|ccccccccc}
	i&5&7&8&4&1&3&6&2   \\\h
	c_{i}(w)&3&4&4&2&0&1&3&0   \\
\end{array}
\]
and so 
$c(w)=0+0+1+2+3+3+4+4=17.$
\end{ex}
\begin{thm}
$c(w)=\dep(w^{-1})$.
\end{thm}
\begin{proof}
Let 
$D_{L}(w)=\{a_{1}, \ds, a_{k}\} (=\Des \,P(w))$ with $a_{1}<\cd<a_{k}$.
Let 
$\a_{1}=a_{1}, 
\a_{j}=a_{j}-a_{j-1} (2\le j\le k), 
\a_{k+1}=n-a_{k}.$
Then 
$\des P(w)=(\a_{1}, \ds, \a_{k+1}) (\te{say} =\a)$
so that 
\[
\dep(w^{-1})=\dep(\des Q(w^{-1}))
=\dep(\des P(w))=\dep\, \al.
\]
On the other hand, 
by definition, the 
sequence $(c_{i}(w))$ looks like 
\[
\begin{array}{c|ccccccccccc}
i&1&\cd&a_{1}&a_{1}+1&\cd&a_{2}&\cd &a_{k}&a_{k}+1&\cd&n   \\\h
c_{i}(w)&0&\cd&0&1&\cd&1&\cd&k-1&k&\cd &k
\end{array}.
\]
Conclude that 
\[
c(w)=\sum c_{i}(w)=
0a_{1}+1(a_{2}-a_{1})+\cd+k(n-a_{k})
=
\dep\, \a.
\]
\end{proof}
\begin{cor}
Charge is Mahonian.
\end{cor}

Thus, the bi-statistic with $(c, \dep)$ directly fits into the framework of RSK correspondence. 
 Recall that for $w\in S_{n}$, the \kyo{inversion} 
statistic 
\[
\el(w)=
\{(i, j)\mid i<j, w(i)>w(j)\}
\]
is also Mahonian, that is, 
$\sum_{w\in S_{n}}t^{\el(w)}=[n]_{t}!$. 
As a consequence, 
\[
\#\{w\in S_{n}\mid \el(w)=k\}=
\#\{w\in S_{n}\mid c(w)=k\}=
\#\{w\in S_{n}\mid \dep\, w=k\}.
\]
Let us observe here a subtle relation 
among these triple Mahonian statistics.
Consider a $(p, q)$-\kyo{bifactorial}
\[
[n]_{p, q}!:=
\sum_{w\in S_{n}}
f_{\l}(p)f_{\l}(q)
=\sum_{w\in S_{n}}p^{c(w)}q^{\dep(w)}.
\]
With $p=1$, it specializes to $[n]_q{!}$  
as the well-known identity says :$\sum_{\l\parn}
f_{\l}f_{\l}(q)=[n]_{q}!$. 
Comparing coefficients of $q^{k}$, we have 
\ytableausetup{mathmode, boxsize=5pt}
\[
\#\{w\in S_{n}\mid \el(w)=k\}=
\sum_{\lparn}f_{\l}\k{[q^{k}](f_{\l}(q))}.
\]
In fact, 
\[
[q^{3}][4]_{p, q}!=
f_{\,\yd{3,1}}(p)+f_{\,\yd{2,1,1}}(p)
=
p+p^{2}+2p^{3}+p^{4}+p^{5}.
\]
\k{\te{missing $f_{\,\yd{2,2}}(p)$}}.
Letting $p=1$
counts 6 permutations in $S_{4}$ with exactly 
3 inversions as 1432, 2341, 2413, 3142, 3214, 4123. 
For a more sophisticated example, 
let us consider 
$[q^{4}][6]_{p,q}!$.
For $\a\models 6$ to have depth 4,
 it must be 24 or 321.
Among all partitions
\[
\te{Par}(6)=\{6, 51, 42, 41^{2}, 3^{2}, 321, 31^{2}, 2^{3}, 2^{2}1^{2}, 21^{4}, 1^{6}\},
\]
there exists a standard Young tableau of 
shape only $51, 42, 41^{2}, 321$ with its descent composition 24 or 321.
\ytableausetup{mathmode, boxsize=15pt}
\[
\begin{ytableau}
	1&2&4&5&6\\
	3\\
\end{ytableau}
\q 
\begin{ytableau}
	1&2&5&6\\
	3&4\\
\end{ytableau}
\q 
\begin{ytableau}
	1&2&3&5\\
	4&6\\
\end{ytableau}
\q
\begin{ytableau}
	1&2&3&5\\
	4\\6
\end{ytableau}
\q
\begin{ytableau}
	1&2&3\\
	4&5\\6
\end{ytableau}
\]
It follows from 
$f_{51}=5, f_{42}=9, f_{41^{2}}=10, f_{321}=16$ that we have 
\[
\#\{w\in S_{6}\mid \el(w)=4\}=
f_{51}+2f_{42}+f_{41^{2}}+f_{321}=49.
\]
Such a counting suggests us a little new problem. 
\begin{prob}
Study a sequence of polynomials 
$\k{[q^{k}][n]_{p, q}!\mid 
0\le k\le \binom{n}{2}
}$. 
For example, when 
are coefficients unimodal?
\end{prob}

\begin{thm}
If $n$ is prime, then 
$[q^{k}][n]_{p, q}!$ has no internal zero.
\end{thm}

\begin{proof}
If $n$ is prime, then all partitions of $n$ are regular. 
Then 
\[
[q^{k}][n]_{p, q}!=
\sum_{\lparn}\k{[q^{k}](f_{\l}(q))}f_{\l}(p)
\]
is a sum of polynomials in $p$ without internal zeros. 
Further, they share a common term $p^{k}$ and all  nonzero coefficients are positive. Thus, 
$[q^{k}][n]_{p, q}!$ has no internal zero.
 \end{proof}
As mentioned in \cite[p.44]{bks}, 
significance of rare internal zeros of fake degree polynomials has been still 
mysterious in algebra and geometry. 
We want to reveal this mystery with perhaps little new combinatorics.


\section{Summary}

In this article, we introduced a new family of polynomials, crystal skeleton polynomials, 
as an auxiliary tool to better understand 
Gessel $F$-expansion of $s_{\l}$, quasi-crystals and crystal skeletons. 
After developing its calculus, we organized four Mahonian statistics through skeleton RSK, fake degree polynomials and a bi-factorial. 
Actually, a part of it was 
a reformulation of previous results in the related literature. 
Nonetheless, we connected and cleared several concepts among little different topics. 
One advantage of using 
descent compositions is 
that we can deal with partitions and 
compositions together in the graded dominance order. 
We wish that our results make one step forward 
to research in related areas.

We plan to establish a theory of 
the skeleton algebra
\[
\mathbf{S}_{n}[x, y]=
\left\{\sum_{w\in S_{n}}a_{w}\,x^{\des w^{-1}}y^{\des w}\Bigg| a_{w}\in\cc\ri\}
\]
and the subalgebra generated by $\Sk_{\l}(x)\Sk_{\l}(y), \lam\vdash n$. Note that we can always recover gradings by 
$x_{i}\mto p^{i-1}x_{i}, y_{i}\mto q^{i-1}y_{i}$. 
A computer experiment suggests that 
$\dim \mathbf{S}_{4}[x, y]=22\ne 4!$ because each of pairs $(1324, 3412), (2143, 4231)$ rarely (but possibly)  shares identical monomials as $x^{22}y^{22}$ and $x^{121}y^{121}$, respectively; see Tables \ref{t3}, \ref{t4}.

\ytableausetup{mathmode, boxsize=14pt}

{\renewcommand{\arraystretch}{1.25}
\begin{table}[h!]
\caption{RS on $S_{4}$ with descent compositions 1}
\label{t3}
\begin{center}
\begin{tabular}{clclclclc}
\toprule
$w$ & $P(w)$ & des\,$P(w)$ & & $Q(w)$ & des\,$Q(w)$ \\
\midrule
\begin{minipage}[t]{0.10\linewidth}1234\end{minipage} & \begin{minipage}[t]{0.18\linewidth}\begin{ytableau}
1 & 2 & 3 & 4
\end{ytableau}\end{minipage} & 4 & & \begin{minipage}[t]{0.18\linewidth}\begin{ytableau}
1 & 2 & 3 & 4
\end{ytableau}\end{minipage} & 4 \\ \addlinespace
\begin{minipage}[t]{0.1\linewidth}1243\end{minipage} & \begin{minipage}[t]{0.18\linewidth}\begin{ytableau}
1 & 2 & 3 \\
4
\end{ytableau}\end{minipage} & 31 & & \begin{minipage}[t]{0.18\linewidth}\begin{ytableau}
1 & 2 & 3 \\
4
\end{ytableau}\end{minipage} & 31 \\ \addlinespace
\begin{minipage}[t]{0.1\linewidth}1324\end{minipage} & \begin{minipage}[t]{0.18\linewidth}\begin{ytableau}
1 & 2 & 4 \\
3
\end{ytableau}\end{minipage} & 22 & & \begin{minipage}[t]{0.18\linewidth}\begin{ytableau}
1 & 2 & 4 \\
3
\end{ytableau}\end{minipage} & 22 \\ \addlinespace
\begin{minipage}[t]{0.1\linewidth}1342\end{minipage} & \begin{minipage}[t]{0.18\linewidth}\begin{ytableau}
1 & 2 & 4 \\
3
\end{ytableau}\end{minipage} & 22 & & \begin{minipage}[t]{0.18\linewidth}\begin{ytableau}
1 & 2 & 3 \\
4
\end{ytableau}\end{minipage} & 31 \\ \addlinespace
\begin{minipage}[t]{0.1\linewidth}1423\end{minipage} & \begin{minipage}[t]{0.18\linewidth}\begin{ytableau}
1 & 2 & 3 \\
4
\end{ytableau}\end{minipage} & 31 & & \begin{minipage}[t]{0.18\linewidth}\begin{ytableau}
1 & 2 & 4 \\
3
\end{ytableau}\end{minipage} & 22 \\ \addlinespace
\begin{minipage}[t]{0.1\linewidth}1432\end{minipage} & \begin{minipage}[t]{0.18\linewidth}\begin{ytableau}
1 & 2 \\
3 \\
4
\end{ytableau}\end{minipage} & 211 & & \begin{minipage}[t]{0.18\linewidth}\begin{ytableau}
1 & 2 \\
3 \\
4
\end{ytableau}\end{minipage} & 211 \\ \addlinespace
\begin{minipage}[t]{0.1\linewidth}2134\end{minipage} & \begin{minipage}[t]{0.18\linewidth}\begin{ytableau}
1 & 3 & 4 \\
2
\end{ytableau}\end{minipage} & 13 & & \begin{minipage}[t]{0.18\linewidth}\begin{ytableau}
1 & 3 & 4 \\
2
\end{ytableau}\end{minipage} & 13 \\ \addlinespace
\begin{minipage}[t]{0.1\linewidth}2143\end{minipage} & \begin{minipage}[t]{0.18\linewidth}\begin{ytableau}
1 & 3 \\
2 & 4
\end{ytableau}\end{minipage} & 121 & & \begin{minipage}[t]{0.18\linewidth}\begin{ytableau}
1 & 3 \\
2 & 4
\end{ytableau}\end{minipage} & 121 \\ \addlinespace
\begin{minipage}[t]{0.1\linewidth}2314\end{minipage} & \begin{minipage}[t]{0.18\linewidth}\begin{ytableau}
1 & 3 & 4 \\
2
\end{ytableau}\end{minipage} & 13 & & \begin{minipage}[t]{0.18\linewidth}\begin{ytableau}
1 & 2 & 4 \\
3
\end{ytableau}\end{minipage} & 22 \\ \addlinespace
\begin{minipage}[t]{0.1\linewidth}2341\end{minipage} & \begin{minipage}[t]{0.18\linewidth}\begin{ytableau}
1 & 3 & 4 \\
2
\end{ytableau}\end{minipage} & 13 & & \begin{minipage}[t]{0.18\linewidth}\begin{ytableau}
1 & 2 & 3 \\
4
\end{ytableau}\end{minipage} & 31 \\ \addlinespace
\begin{minipage}[t]{0.1\linewidth}2413\end{minipage} & \begin{minipage}[t]{0.18\linewidth}\begin{ytableau}
1 & 3 \\
2 & 4
\end{ytableau}\end{minipage} & 121 & & \begin{minipage}[t]{0.18\linewidth}\begin{ytableau}
1 & 2 \\
3 & 4
\end{ytableau}\end{minipage} & 22 \\ \addlinespace
\begin{minipage}[t]{0.1\linewidth}2431\end{minipage} & \begin{minipage}[t]{0.18\linewidth}\begin{ytableau}
1 & 3 \\
2 \\
4
\end{ytableau}\end{minipage} & 121 & & \begin{minipage}[t]{0.18\linewidth}\begin{ytableau}
1 & 2 \\
3 \\
4
\end{ytableau}\end{minipage} & 211 \\
\addlinespace
\end{tabular}
\end{center}
\end{table}}

{\renewcommand{\arraystretch}{1.25}
\begin{table}[h!]
\caption{RS on $S_{4}$ with descent compositions 2}
\label{t4}
\begin{center}
\begin{tabular}{clclclclc}
\toprule
$w$ & $P(w)$ & des\,$P(w)$ & & $Q(w)$ & des\,$Q(w)$ \\
\midrule
\begin{minipage}[t]{0.1\linewidth}3124\end{minipage} & \begin{minipage}[t]{0.18\linewidth}\begin{ytableau}
1 & 2 & 4 \\
3
\end{ytableau}\end{minipage} & 22 & & \begin{minipage}[t]{0.18\linewidth}\begin{ytableau}
1 & 3 & 4 \\
2
\end{ytableau}\end{minipage} & 13 \\ \addlinespace
\begin{minipage}[t]{0.1\linewidth}3142\end{minipage} & \begin{minipage}[t]{0.18\linewidth}\begin{ytableau}
1 & 2 \\
3 & 4
\end{ytableau}\end{minipage} & 22 & & \begin{minipage}[t]{0.18\linewidth}\begin{ytableau}
1 & 3 \\
2 & 4
\end{ytableau}\end{minipage} & 121 \\ \addlinespace
\begin{minipage}[t]{0.1\linewidth}3214\end{minipage} & \begin{minipage}[t]{0.18\linewidth}\begin{ytableau}
1 & 4 \\
2 \\
3
\end{ytableau}\end{minipage} & 112 & & \begin{minipage}[t]{0.18\linewidth}\begin{ytableau}
1 & 4 \\
2 \\
3
\end{ytableau}\end{minipage} & 112 \\ \addlinespace
\begin{minipage}[t]{0.1\linewidth}3241\end{minipage} & \begin{minipage}[t]{0.18\linewidth}\begin{ytableau}
1 & 4 \\
2 \\
3
\end{ytableau}\end{minipage} & 112 & & \begin{minipage}[t]{0.18\linewidth}\begin{ytableau}
1 & 3 \\
2 \\
4
\end{ytableau}\end{minipage} & 121 \\ \addlinespace
\begin{minipage}[t]{0.1\linewidth}3412\end{minipage} & \begin{minipage}[t]{0.18\linewidth}\begin{ytableau}
1 & 2 \\
3 & 4
\end{ytableau}\end{minipage} & 22 & & \begin{minipage}[t]{0.18\linewidth}\begin{ytableau}
1 & 2 \\
3 & 4
\end{ytableau}\end{minipage} & 22 \\ \addlinespace
\begin{minipage}[t]{0.1\linewidth}3421\end{minipage} & \begin{minipage}[t]{0.18\linewidth}\begin{ytableau}
1 & 4 \\
2 \\
3
\end{ytableau}\end{minipage} & 112 & & \begin{minipage}[t]{0.18\linewidth}\begin{ytableau}
1 & 2 \\
3 \\
4
\end{ytableau}\end{minipage} & 211 \\ \addlinespace
\begin{minipage}[t]{0.1\linewidth}4123\end{minipage} & \begin{minipage}[t]{0.18\linewidth}\begin{ytableau}
1 & 2 & 3 \\
4
\end{ytableau}\end{minipage} & 31 & & \begin{minipage}[t]{0.18\linewidth}\begin{ytableau}
1 & 3 & 4 \\
2
\end{ytableau}\end{minipage} & 13 \\ \addlinespace
\begin{minipage}[t]{0.1\linewidth}4132\end{minipage} & \begin{minipage}[t]{0.18\linewidth}\begin{ytableau}
1 & 2 \\
3 \\
4
\end{ytableau}\end{minipage} & 211 & & \begin{minipage}[t]{0.18\linewidth}\begin{ytableau}
1 & 3 \\
2 \\
4
\end{ytableau}\end{minipage} & 121 \\ \addlinespace
\begin{minipage}[t]{0.1\linewidth}4213\end{minipage} & \begin{minipage}[t]{0.18\linewidth}\begin{ytableau}
1 & 3 \\
2 \\
4
\end{ytableau}\end{minipage} & 121 & & \begin{minipage}[t]{0.18\linewidth}\begin{ytableau}
1 & 4 \\
2 \\
3
\end{ytableau}\end{minipage} & 112 \\ \addlinespace
\begin{minipage}[t]{0.1\linewidth}4231\end{minipage} & \begin{minipage}[t]{0.18\linewidth}\begin{ytableau}
1 & 3 \\
2 \\
4
\end{ytableau}\end{minipage} & 121 & & \begin{minipage}[t]{0.18\linewidth}\begin{ytableau}
1 & 3 \\
2 \\
4
\end{ytableau}\end{minipage} & 121 \\ \addlinespace
\begin{minipage}[t]{0.1\linewidth}4312\end{minipage} & \begin{minipage}[t]{0.18\linewidth}\begin{ytableau}
1 & 2 \\
3 \\
4
\end{ytableau}\end{minipage} & 211 & & \begin{minipage}[t]{0.18\linewidth}\begin{ytableau}
1 & 4 \\
2 \\
3
\end{ytableau}\end{minipage} & 112 \\ \addlinespace
\begin{minipage}[t]{0.1\linewidth}4321\end{minipage} & \begin{minipage}[t]{0.18\linewidth}\begin{ytableau}
1 \\
2 \\
3 \\
4
\end{ytableau}\end{minipage} & 1111 & & \begin{minipage}[t]{0.18\linewidth}\begin{ytableau}
1 \\
2 \\
3 \\
4
\end{ytableau}\end{minipage} & 1111 \\ \addlinespace
\end{tabular}
\end{center}
\end{table}}


%





\end{document}